\documentclass[twoside,10pt]{article}
\usepackage[spanish, mexico]{babel}
\usepackage[latin1]{inputenc}
\usepackage{multicol}
\usepackage{graphicx}
\usepackage{amssymb}
\usepackage{amsmath,amsfonts}
\usepackage{enumerate}
\usepackage{array}
\usepackage{ textcomp } 
\usepackage{mathtools}
\usepackage{color}
\usepackage[list=true, font=large, labelfont=bf, 
labelformat=brace, position=top]{subcaption}
\usepackage[]{todonotes}
\usepackage{tikz}
\usetikzlibrary{knots}

\input xy
\xyoption{all}

\newtheorem{theorem}{Theorem}[section]
\newtheorem{lemma}[theorem]{Lemma}

\newtheorem{corollary}[theorem]{Corollary}

\newtheorem{definition}[theorem]{Definition}
\newtheorem{example}[theorem]{Example}

\newtheorem{remark}[theorem]{Remark}

\makeatletter
\def\blfootnote{\gdef\@thefnmark{}\@footnotetext}
\makeatother

\usepackage{microtype}
\usepackage{xcolor}
\usepackage{lmodern}
\usepackage{graphicx}
\usepackage{listings}
\usepackage{paralist}
\usepackage{setspace}
\usepackage{calc}
\usepackage{siunitx}
\usepackage{multirow}
\usepackage{threeparttable}
\usepackage{booktabs}
\usepackage{float}
\usepackage{framed}
\usepackage{titling}
\usepackage{esvect}
\usepackage{tikz}
\usepackage{pgfplots}
\usetikzlibrary{intersections}
\usetikzlibrary{shapes,arrows}
\usetikzlibrary{calc}
\usetikzlibrary{shapes.geometric}
\usetikzlibrary{shapes.arrows}
\usetikzlibrary{arrows.meta}
\usetikzlibrary{decorations.markings}
\usetikzlibrary{fit}
\usetikzlibrary{hobby}

\pgfplotsset{compat=1.12}
\newcommand{\C}{\mathbb{C}}
\newcommand{\R}{\mathbb{R}}
\newcommand{\Z}{\mathbb{Z}}

\numberwithin{equation}{section}

\begin{document}

\thispagestyle{plain} 


\bigskip

\begin{center}

{\bf {\Large Complete hyperbolic structures in the complement of the Borromean rings.}}

\vspace{0.2cm}

{\large {$^a$ Angel Cano, $^b$ Juan Francisco  Estrada}}

\vspace{0.2cm}

$^{a}$ Unidad Cuernavaca del Instituto de Matemáticas de la UNAM.

$^{b}$ Facultad de Ciencias Fisico Matem\'aticas. Benémerita Universidad Autónoma de Puebla.

\medskip

$^a$angelcano@im.unam.mx, $^b$ festrada@fcfm.buap.mx

\end{center}

\renewcommand{\abstractname}{Abstract}
\begin{abstract}
	In this note, we show the fundamental group of the complement of the Borromean rings in $\Bbb{S}^3$ has exactly two representations in ${\rm PSL}(2,\Bbb{C})$ which are faithful, discrete and send meridians into parabolic elements. Using this result we are able to show that Borromean rings admit exactly one hyperbolic structure.
	
\end{abstract}


$\underline{ \ \ \ \ \ \ \ \ \ \ \  \ \ \ \ \ \ \ \ \ \ \ \ \ \ \ \
\ }$

{\footnotesize $Keywords$ $and$ $phrases:$ Quandles,  Representation Variety, Links.}

{\footnotesize 2010 $Mathematics$ $Subject$ $Classification$ Primary 37F99, 32Q, 32M Secondary 30F40, 20H10, 57M60, 53C.}

$\overline{ \ \ \ \ \ \ \ \ \ \ \  \ \ \ \ \ \ \ \ \ \ \ \ \ \ \ \ \
}$

\section*{Introduction}We can say that Ahlfors'  finiteness theorem, see \cite{ahl1, ahl2},   and Sullivan's finiteness theorem for cusps, see \cite{sul1}, are the major pillars in   Kleinian groups,  with deep implications in dynamics, topology and geometry, see for example \cite{Hubbard,sul2}. For several years after the publication of this results, mathematicians around the world worked on new proofs with the aim to generalize this theorems for higher dimensional conformal Kleinian groups, see \cite{kul}, despite the known evidence in those years M. Kapovich and L.  Potyagalio constructed an example that shows that this kind of theorems cannot hold in the higher dimensional conformal setting, see \cite{kapovich}.  The  construction Kapovich and  Potyagalio is based in the existence of a complete hyperbolic structure on the complement of the Borromean rings, more precisely Kapovich-Potyagalio utilizes a discrete subgroup of Möbius transformation whose quotient in $\Bbb{H}^3$  is the complement of the Borromean rings in $\Bbb{S}^3$.  The main  purpose is to construct all the  reprsentations of the fundamental group of the complent of the Borromean ring into $Mob(\hat{C})$ which sends meridians into parabolic elements, more precisely in    this article we show:

\begin{theorem} \label{t:main1}
	There are exactly two representations, namely  $\rho_1,\rho_2$,  from the complement of the Borromean rings into ${\rm PSL(3,\Bbb{C})}$ which are faithful, have a discrete image and sends meridians into parabolic elements.  Moreover, we have $\rho_1$ and $\rho_2$ are complex conjugate each other, {\it i. e.} $\rho_1=\overline{\rho_2}$.
\end{theorem}

\begin{theorem} \label{t:ext}
	The complement of the Borromean in $\Bbb{S}^3$ (as well as its mirror image) has exactly one complete real hyperbolic structure. Moreover if $\rho_ 1$ , $\rho_2$ are the respective holonomies of the complete real hyperbolic structures on the Borromean links and its mirror image, then $\rho_1=\overline{\rho_2}$. 
\end{theorem}

The paper is organized as follows: Section \ref{s:pre} reviews some well-known facts on Knots and Quandles, see \cite{Joyce,rolf} for detailed information. In Section \ref{s:col}, we compute a representation for the fundamental quandle of the Borromean rings and as a consequence, we are able to compute all ${\rm PSL(2,\Bbb{C})}$-colorations for the rings.  Finally in section \ref{s:main}  we prove the main results of this article.  The results exposed here well know, see \cite{thurston} for a beautiful sketch and  \cite{ucan2} for a  detailed proof, and the procedure is not new,  as far as we know it was introduced by Inoue and Kabaya in \cite{IK}, but here is presented in a systematic manner so it can be used for future reference.  The article is itself contained so we hope this can be read by a wide audience.

\section{Preliminaries} \label{s:pre}

\subsection{Knots}
In this subsection we review the facts used trough the text about links and knots, there are tons of excellent books on this subject,  we find useful and very well writing the following books  \cite{BZ, man, rolf}.    A {\bf knot} is an embedding of a circle $\Bbb{S}^1$ into the 3-sphere, $\Bbb{S}^3$,  more generally a  {\bf link} is a collection of knots which do not intersect.  A knot is called {\bf tame} if it is isotopic to a
{\bf polygonal knot}, {\it i. e. } a knot in $\Bbb{R}^3$ which is the union of a finite set of line segments (the endpoints of such segment lines are called the vertices). Similarly, we will say that a link is tame if it is isotopic to a link where each connected component is a polygonal knot, trough this article all the links are assumed to be tame. An {\bf oriented link} is a link where each connected component has an orientation. For example the {\bf trivial knot} is given by the following embedding $\mathbb{S}^1\rightarrow \mathbb{R}^3$ given by $e^{2 \pi i \theta}\mapsto (Cos \theta,Sen \theta, 0).$

A useful way to visualize and manipulate links is to project the link onto a plane, more precisely  if $L\subset \R^3$ is a polygonal link, $P\subset \R^3 $  a plane and    $p:\R^3\rightarrow P$  is the orthogonal projection, the projection  $p$ is said to be {\bf regular}  
if: 
\begin{enumerate}
	\item   The set of points  $x \in  p(L)$   whose preimage $p^{-1}(x )$ under the projection contains more than one points (crossing points) is finite and the crossing points are double points,
	\item  no vertex of $L$ is mapped onto a double point.
\end{enumerate}
It is possible to show given a polygonal link we can find an isotopic link with a regular projection, see \cite{BZ}.  A {\bf link diagram} is the regular projection of a link  to the plane with broken
lines indicating where one part of the link under crosses the other part,
a connected component of a link diagram will be called an {\bf arc} inherit the orientation also it is trivial that a    link can be reconstructed from its diagram, see \cite{BZ}. 

In order to determine whether two projections come from the same link,  Reidemeister introduced three transformations, called {\bf  Reidemeister moves}, that change the projection but don't change the knot: the first one adds or removes a link,   the second pass or underpass strand of the knot over another and the third passes a strand over a crossing of another, see the figure below.

\begin{figure} [!htb] 
	\begin{center}
		\includegraphics[width=60mm]{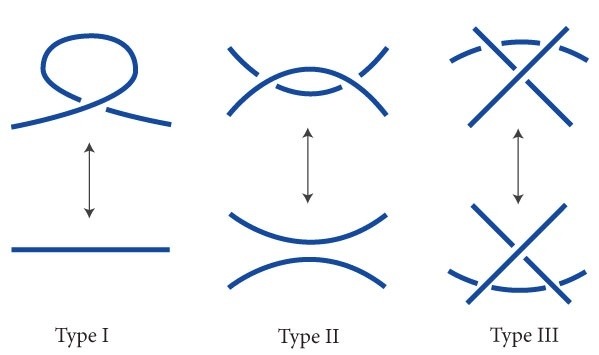}
	\end{center}
\end{figure}

\begin{theorem}[Reidemaster, see \cite{man}]
	Two links are equivalent if and only if all their diagrams are connected by a finite sequence of Reidemeister moves.
\end{theorem}

Observe that in order to prove the invariance for some function on knots, one should check its invariance under Reidemeister moves. 

The following result due to Wirtinger provides us a way to compute a  presentation for the fundamental group of a link complement.   

\begin{theorem}[Wirtinger, see \cite{rolf} ]\label{t:wir}
	Let $L$ be a  link  expressed by a link diagram, and let $A$ be the set of arcs and $B$ the set of crossings. Let $W$ be the free group with generating set $A$, and let $N$ be the subgroup of $A$ generated by the elements $r(b)$ for each $b \in  B$, with $r$ refined as follows:
	\[
	r(b)=
	\left \{
	\begin{array}{ll}
		(u(b)+1)o(b)u(b)^{-1}o(b)^{-1}    & \textrm{ if  } b \textrm{ is right-handed}\\
		o(b)(u(b) + 1)o(b)^{-1}u(b)^{-1} & \textrm{ if  } b \textrm{ is left-handed}\\
	\end{array}
	\right. 
	\]
	here $u(b)+1, u(b),o(b)$, right  and left handed  are as in the picture bellow. Then $W/N$ is isomorphic to the fundamental group of the complement of $L$ in $\Bbb{S}^3$.
\end{theorem}

\begin{figure} [!htb] 
	\begin{center}
		\includegraphics[width=80mm]{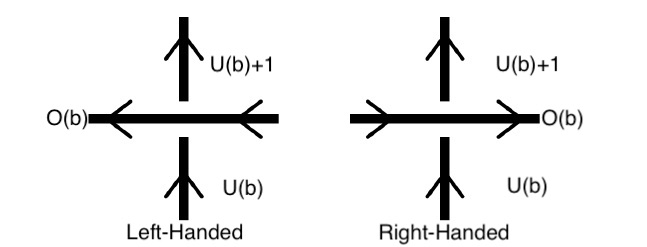}
	\end{center}
\end{figure}
We must point out, to the novice in Knot theory, that decide whether two knots given knots are equivalent is not an easy task.  And maybe one of the simplest criteria to decide is to knots are candidates to be the same is the so-called {\bf tricolorability}, a  knot is {\bf tricolorable} if each arc  of the knot diagram can be colored one of three colors, subject to the following rules:
\begin{enumerate}
	\item  All colors must be used, and
	\item  At each crossing, the three incident strands are either all the same color or all different colors.
\end{enumerate}

Trough the previous notion we can ensure the trefoil knot, see figure below, is not equivalent to the trivial knot.
\begin{figure} [!htb] 
	\begin{center}
		\includegraphics[width=40mm]{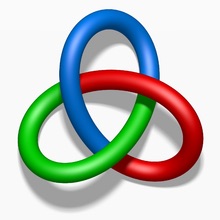}
	\end{center}
\end{figure}

Finally, a {\bf longitude } is a loop that runs from a  basepoint $x$ along a path to a point y on the boundary of a tubular neighborhood of the link, then follows along the tube, making one full lap to return to y, then returns to x via the path. A {\bf meridian} is a loop that runs from x to y, then circles around the tube, returns to y, then returns to x. (The property of being a longitude or meridian is well-defined because the tubular neighborhoods of a tame knot are all ambiently isotopic).

\subsection{Quandles}
The Quandles were introduced by Joyce in \cite{Joyce} as a way to respond a question of Fox about peripheral groups of knot groups, and quite recently they have gained great attention for their application in several areas of the math. Most of literature about Quandles 
is dispersed in articles and as a  chapter in books on Knot theory, the interested reader can see  \cite{nosaka} for interesting applications of quandles and an excellent review of the basic theory.

\begin{definition}
	A quandle   is a set $X$ with a binary operation $*:X\times X\rightarrow X$ satisfying:
	\begin{enumerate}
		\item  \label{i:q1} For all $x\in X, \, x* x=x$.
		
		\item \label{i:q2} For all $y\in X$, the map  $\beta_y: X\rightarrow X$ defined by $\beta_y(x)=x* y$ is invertible.
		
		\item \label{i:q3} For all $x,y,z \in X, \, (x * y)* z =(x*z)*(y * z) $.
		
	\end{enumerate}
	
\end{definition}

In the following we will write  $x *^{-1} y$ for $\beta^{-1}_y(x)$.
The quandle axioms (\ref{i:q2}) and (\ref{i:q3}) are equivalent to saying that $$Inn(X)=\{\beta_y:X\rightarrow X: y\in X \}$$ is a group, so we can define a map   $inn: X \rightarrow Inn(X)$ by $y\mapsto \beta^{-1}_y$.
A map $f : X \rightarrow  Y$ between two quandles  X, Y is called a {\bf quandle
	morphism} if $f(a * b) = f(a) * f(b)$ for any $a, b \in  X$.

\begin{example}[Trivial quandle]
	Any non-empty set $X$ with the trivial operation $x* y=x$ is a quandle.
\end{example}

\begin{example}[n-fold conjugation quandle]
	Let $G$ be any group, and $H\subset G$  a set  closet under conjugations, then $H$ is a quandle  under the $n$-fold conjugation {\it i. e. }  $x* y=y^{-n} xy^{n}$, in the following we will denote this quandle by $Conj_n(H)$. 
\end{example}

\begin{example}[Parabolic quandle] \label{e:par}
	The set $(\Bbb{C}^2 \setminus \{0\})/\pm $ with the binary operation  
	\[
	[x_1,y_1]*[x_2,y_2] = [x_1+x_1x_2y_2-y_1x_2^2, x_1 y_2^2 +y_1- x_2y_2 y_1]
	\]
	is a quandle, see \cite{IK}, this quandle will be denote by $Par$. 
\end{example}

It is straightforward to check that the parabolic quandle is isomorphic to the 1-fold conjugation quandle of parabolic transformations of $PSL(2,\Bbb{C})$.

\begin{example}[$n$-Dihedral quandle]
	Let $n$ be a positive integer,  for $ i, j \in \{0, 1, . . . , n - 1\}$,
	define $i *  j= 2j - i$ mod $n$,  this defines  a quandle. This quandle can be viewed also as the 1-fold conjugation quandle   of reflections in the dihedral group.
\end{example}

\begin{example}[Eisermann quandle]
	The disjoint union $\mathcal{Q}_{m,n} = \Z_m \sqcup \Z_n$ becomes a quandle with
	$a * b = a$ for $a,b \in \Z_m$ or $a,b \in \Z_n$, and $a *b = a + 1$ otherwise. 
\end{example}

\begin{example}[Free quandle]
	Given an index set $I$ the Free quandle associated with $I$,  in symbols $Q^{F}_I$, is the set of conjugacy classes of  $I$ in the free group generated by $I$.
\end{example}

As it happens for the case of groups free quandle is a free object, {\it i. e. } if $I$ is an index set and $\alpha: I \rightarrow X$ is a function, then there exist a unique quandle morphism $f: Q^F_I\rightarrow X$  such that $f(y)=\alpha(y)$ for any $y\in I$.

\begin{example}[Fundamental Quandle-Algebraic Version]
	Let  $D$ be a  diagram of an oriented link  $K$.   Denote the set of arcs of $D$ by $A_D$ and the set of croossings by $C_D$,  let $Q^F_{A_D}$ be the free quandle assosiated to $A_D$, we introduce the equivalence relation $\sim $ for $Q^F_{A_D}$ as follows:
	\begin{enumerate}
		\item $x\sim x$ for every $x\in Q^F_{A_D}$.
		\item If  $p\in C_D$ is incident to two undercrossing arcs $a_p$ and $c_p$ and an overcrossing arc $b_p$,  then  $a_p * b_p \sim  c_p$ and $c_p\sim a_p * b_p $, where $a$ is the arc lying on the left hand with respect to $b$ and $c$ is the arc lying on the right hand with respect to $b$.
		\item
		For every  $w_1,w_2 \in  Q^F_{A_D}$  we have  $(w_1* (a_p*b_p))*w_2\sim (w_1*c_p)*w_2$ and $(w_1*c_p)*w_2\sim (w_1* (a_p*b_p))*w_2$ for all $p\in C_D$. 
		
	\end{enumerate}   The set of equivalence classes is denoted by $\Gamma\langle A_D|R_D\rangle$ and 
	it is easy to check that it is a quandle with respect to the operation $[W_1]*[W_2]=[W_1*W_2]$, this will called the fundamental quandle of $K$. 
\end{example}

As the reader has noted quandles work in a similar way as groups.

\begin{example}[Fundamental Quandle (Geometric Version)]
	Let $N(K)$ a regular neighborhood of $K$, then the fundamental
	quandle of $K$ is the set of homotopy classes of paths in $\Bbb{S}^3\setminus  N(K)$ from a base point to $N(K)$ such
	that the initial point stays fixed at the base point while the terminal point is free to wander on
	$N(K)$, the quandle operation is then given by setting $x * y$ to the homotopy class of the path given
	by first following $y$, then going around a canonical meridian on $N(K)$ linking $K$ once, then going
	backward along $y$, then following $x$. 
\end{example}

\begin{theorem}
	The fundamental quandle of a link given by the geometric version agree with the algebraic one.  
\end{theorem}
\begin{definition}
	Given a quandle $X$ we define its adjoint group $Adj(X)$ to be
	the quotient group of the free group  generated by the set $X$ modulo the relations
	induced by the quandle operation,  $R=\{ a*b =b^{-1}ab :a,b\in X\}$. By construction we
	obtain a  map  $adj: X\rightarrow Adj(X)$
	with  $adj(a*b)=adj(b)^{-1}adj(a)adj(b)$.
\end{definition}

The group $Adj(X)$ can be interpreted as the enveloping group of $X$. However,
the map $adj$  could or could not be injective, see example  \ref{e:injective} and \ref{e:noinjective}  below.

\begin{example}\label{e:injective} Given a group $G$ we have 
	$Adj(Conj_1(G))=G.$
\end{example}
\begin{example}[See \cite{Eisermann2}]\label{e:noinjective} 
	For  $m,n \in  \Bbb{N}$ the  group $Adj(\mathcal{Q}_{m,n})$ is isomorphic
	to the quotient $Heis(3, \Bbb{Z}) / \langle z^l \rangle$, where  
	$l= gcd(m,n)$ and 
	
	\[
	z=
	\left (
	\begin{array}{lll}
		1 & a &b\\
		0 & 1 & c\\
		0 & 0 & 1
	\end{array}
	\right )
	\]
	
	In particular, $adj: \mathcal{Q}_{m,n}\rightarrow  Adj(\mathcal{Q}_{1m,n})$ is injective if and only if $m = n$.
\end{example}
In general it is hard to determine explicitely  the adjoint group  of a quandle  but in some cases it is possible to determine it, for example 
in the case of knots the adjoint group, as expected from Wirtinger theorem see \ref{t:wir},  the associated group is  simply the fundamental group, more precisely:

\begin{theorem}[See \cite{Joyce}]
	Let $\mathcal{L}\subset \Bbb{S}^3$ be a link, with diagram $D$ and set of  arcs $A_D$, thus $Adj(\Gamma\langle A_D|R_D\rangle )$ is isomorphic to $\pi_1( \Bbb{S}^3\setminus \mathcal{L})$.
\end{theorem}

The previous theorem and the geometric interpretation of the fundamental quandle enable us to construct,  in a natural way, the map  $adj: \Gamma\langle A_D|R_D\rangle \rightarrow  \pi_1(\Bbb{S}^3\setminus \mathcal{L})$ from the link quandle to the link group as follows: to an element $\gamma$ in the quandle associate the loop $\gamma m \gamma^{-1}$, where $m$ is the meridian at the point $x$. Observe that this shows thas $adj$ send arcs into meridians.

The adjoint group plays a relevant a role in this article, let us summarize some of its properties, here we include a proof because it is not availble in the litterature and in this article will be necessary to understand how to construct some of the maps appearing in the theorem.

\begin{theorem}[See \cite{Eisermann2}] \label{p:biy} Let $X,Y$ be two quandles and $G$ a group,  then:
	\begin{enumerate}
		\item   \label{ad1} {\bf Universal property.} If   $\phi : X\rightarrow G$  satisfy $ \phi (a * b) = \phi(b)^{-1}\phi(a)  \phi(b)$ for all $a,b \in  X$, then  there exists a
		unique group morphism $h: Adj(X) \rightarrow  G$ such that $\phi  = h \circ  adj$.
		
		\item  {\bf Functoriality.} Every quandle morphism $\psi: X \rightarrow Y$
		induces a unique group morphism $Adj(\psi ): Adj(X) \rightarrow Adj(Y)$ so the following diagram commutes:
		
		\begin{equation} 
			\xymatrix{
				X \ar[r]^{\psi} \ar[d]^{adj_X}         & Y \ar[d]^{adj_Y}\\
				Adj(X) \ar[r]^{Adj(\psi)}  & Adj(Y)
			} 
		\end{equation} $Adj(\psi)adj_X = adj_Y\circ \psi$.

		\item  {\bf Adjointness}. We
		have a natural bijection $$Hom_{Qnd}(Q,Conj_1(G)) \equiv Hom_{Grp}(Adj(Q),G).$$

		\item {\bf Adjoint action}. The map  $inn: X\rightarrow Inn(X)$ which is given by $inn(a)(b)=b*^{-1}a$ induces a unique
		group homomorphism $\rho_X  : Adj(X) \rightarrow Inn(X)$ such that $inn = \rho_X  \circ adj$. In this way the adjoint
		group $Adj(X)$ defines a right action on  $X$  by $xg=\rho_X(g)(x)$,  more explicitly if $g=x_1^{\epsilon_1}x_2^{\epsilon_2}\cdots x_n^{\epsilon_n}\in Adj(X)$  where   $n\geq 0$, $x_i \in X$, and $\epsilon\in\{\pm 1\}$, thus
		\[
		x*g=(\ldots((x*^{\epsilon_1}x_1)*^{\epsilon_2}x_2)\ldots)*^{\epsilon_n}x_n.
		\]

		\item {\bf Equivariance}.  Let   $\phi : X \rightarrow  Y$  be  a quandle morphism 
		and $adj(\phi): adj(X)\rightarrow adj(Y)$   the induced group  morphism, then we can  define a right action of  $Adj(X)$  on  $Y$
		by $y*g=y*adj(\phi)(g)$, making $\phi$  equivariant under the right action of $Adj(Y)$.
		
	\end{enumerate}
\end{theorem}

{\it Proof.- }
Let us prove the universal property. Let  $\phi : X\rightarrow G$  be  map satisfying $\phi(a*b)=\phi(b)^{-1}\phi (a)\phi(b)$, by the universal property of the free group  there is a group morphism $\hat \phi $  such that $\hat \phi \circ i=\phi $, where $i: X\rightarrow Free(X)$ is the canonical inclusion.  Now observe that for every $a,b\in X$ we have:
$$\hat \phi ((a*b)b^{-1}a^{-1}b)=\hat \phi (a*b)\hat\phi (b)^{-1}\hat \phi( a)^{-1}\hat \phi( b)=  \phi (a*b)\phi (b)^{-1} \phi( a)^{-1} \phi( b)=e$$
Thus $R=\langle  (a*b)b^{-1}  a^{-1}  b\vert a,b\in x \rangle\subset Ker(\hat \phi) $, thus  by the universal property of the quotient group there is a group morphish $h:Free(X)/R=Adj(X)\rightarrow G$ such that  $h \circ q=\hat \phi $, here $h: Free (X)\rightarrow Free(X)/R $ denotes the quotient map.  Observe  that $\phi=\hat \phi \circ i=  h\circ q\circ i=h\circ adj  $, which shows the existence of the map.  Now, let us assume there is another group morphism $\tilde d: Adj(x)\rightarrow G$ so that $\tilde d\circ  adj=\phi$,   observe that for every $a\in X$ we have: 
$$\tilde d\circ q(x)=\tilde d\circ  q \circ i(a)=\tilde d\circ adj(x)= \phi(x)=d\circ adj(x)=d\circ q \circ i(x)=d\circ q(x)=\hat \phi (x)$$ so $\hat \phi=\tilde d\circ q$,  by the universal property of the quotient group we have $d=\tilde d$ which shows the uniqnees.

The functoriality will be easy from the previous proposition. Observe that if $\psi: X\rightarrow Y$ is a quandle morphism thus $adj\circ \psi: X\rightarrow Adj(Y)$ is a map that satisfies: $$adj\circ \psi (a*b)=adj(\psi (a)*\psi (b))=adj\circ \psi  (b)^{-1}adj\circ\psi  (a)adj\circ\psi  (b)$$
so by the universal property of the adjoint group there is a unique group morphism $Adj(\psi)$ such that $Adj(\psi)\circ adj=adj \circ \psi $, which concludes this part or the proposition.

In order to show adjointness let us define $Adj: Hom_{Qnd}(Q,Conj_1(G)) \rightarrow  Hom_{Grp}(Adj(Q),G)$ and $Conj:   Hom_{Grp}(Adj(Q),G) \rightarrow Hom_{Qnd}(Q,Conj_1(G)) $ as follows:  if $\psi: Q\rightarrow Conj_1(G)$ is a quandle morphism, then $Adj(\psi): Adj(Q)\rightarrow G$  is the group morphism given by the functoriality property, now if $\phi : Adj(Q)\rightarrow G$  is a group morphism  we define  $Conj(\phi)=\phi \circ adj: Q\rightarrow G$ which is quandle morphism since 
$$
\begin{array}{ll}
	Conj(\phi)(a*b)&=\phi \circ adj(a*b)\\
	&=\phi (adj(b)^{-1}abj(a)adj(b))\\
	&=\phi (adj(b)^{-1})\phi (adj(a))\phi(adj(b))\\
	&=Conj(\phi)(b)^{-1}Conj(\phi)(a)Conj(\phi)(b)
\end{array}
$$ 

Now let $\psi\in Hom_{Qnd}(Q,Conj_1(G))$, then $Conj(Adj(\psi))=Adj(\psi)\circ adj=adj\circ\psi=\psi $. Finally, if  $\phi\in Hom_{Grp}(Adj(Q),G)$, then $Adj(Conj(
\phi))=Adj(\phi\circ adj )=\phi $, which shows this part of the  proposition. 

Now, by applying the universal property and the following straightforward computations we are to show the adjoint action property.
\[
\begin{array}{l}
	inn(a)*inn(b)(c)=
	(c*^{-1}b)*^{-1}a)*b
	= c*^{-1}(a*b)=inn(a*b)(c)
\end{array}
\]

The last proposition is straightforward so we will omit its proof here.
$\square$

As a corollaries we have
\begin{corollary} \label{c:mod}
	Let $X$  be quandle and $G$ be a group. Let us consider the left action of $G$ by conjugation on  $Hom_{Grp}(Adj(X),G)$ and 
	$Hom_{Qnd}(X,Conj_1(G))$, thus we have  
	
	\[
	Hom_{Qnd}(X,Conj_1(G)) /G\equiv Hom_{Grp}(Adj(X),G)/G. 
	\]   
\end{corollary} 
{\it Proof.-}  Let us show $Conj:   Hom_{Grp}(Adj(X),G) \rightarrow Hom_{Qnd}(X,Conj(G)) $, is $G$-equivariant. Let $g\in G$ and   $\phi \in Hom_{Grp}( Adj(X),G)$, then  
\[
Conj(g\phi)(c)=g\phi \circ adj(c)=g^{-1}\phi (adj(c))g=Conj(\phi)(c)*g,
\] 
this  concludes the proof.
$\square$

\begin{definition}
	An {\bf arc coloring} of $\Gamma\langle A_D|R_D\rangle$ in $X$ is a quandle morphism  $A : \Gamma\langle A_D\vert R_D \rangle  \rightarrow X$.
\end{definition}
\begin{example}
	If $X$ is any quandle and $x\in X$, then the constant function $y\mapsto x$ is an arc coloring.
\end{example}

\begin{remark}
	\begin{enumerate}
		\item[ ] 
		\item 
		Tricolarability of links is equivalent to provide surjective arc coloring from the fundamental quandle to the 3-dihedral quandle. 
		\item 
		Given a link   $L\subset \Bbb{S}^3$,  with projection  $D$ and a quandle  $X$ quandle, in order to produce an arc coloring of $D$ in $X$,    we may  choose an image $f(x_k)\in X$ for every arc $x_k$, then  once we define  $f(x_k)$ and $f(x_j)$, we simply put $f(x_k *x_j)=f(x_k) * f(x_j)$. However we must be cautious, since  potential problem arises at the crossings, so if the arcs satisfy the relation $x*y = z$ at a crossing,  we should  choose $f(x), f(y), f (z)$ so that
		$ f (x) *f (y) =f (z)$. 
	\end{enumerate} 
\end{remark}

\begin{corollary} \label{c:correp}
	Let $\mathcal{L}$ be an oriented link in $\Bbb{S}^3$, $D$ a diagram of $\mathcal{L}$, and $A$ an arc coloring of $D$ with respect to $Par$, then   we have a representation $\rho : \pi_1(\Bbb{S}^3\setminus \mathcal{L}) \rightarrow PSL(2,\C),$
	derived from $A$ that  sends meridians to parabolic elements.  Conversely, given   $\rho  : \pi_1(\Bbb{S}^3\setminus \mathcal{L})\rightarrow PSL(2,\Bbb{C})$  a representation  sending meridians into parabolic elements, then  we can define an arc coloring $A$ of $D$ respect to $Par$ such that the induced representation into $PSL(2,\Bbb{C})$ is $\rho$.
\end{corollary}

{\it Proof.-}
Let $A: \Gamma\langle D\vert R\rangle\rightarrow Par $ be a coloration since $Par$ is isomorphic to the 1-fold conjugation quandle of parabolic transformations of $PSL(2,\Bbb{C})$, we have $A$ can be thinked as a coloration of the fundamental knot into  $PSL(2,\Bbb{C})$ that sends arcs into parabolic elements. By the universal property  there is a unique  representation   $\rho:\Pi(\Bbb{S}^3\setminus \mathcal{L}) \rightarrow PSL(2,\Bbb{C})$ such that $A=\rho\circ adj$. Since $adj$ sends arcs into meridians and $A$ sends arcs we conclude $\rho$ sends meridians into parabolics. Now the rest of the proof is trivial.
$\square$

\section{ Borromean rings   ${\rm PSL}(2,\Bbb{C})$-parabolic colorations} \label{s:col}
The so-called Borromean rings are three circles that are linked in such a way that when we remove one of the circles the other two are not linked, see figure below. 

\begin{figure} [!htb] \label{fig:boat1} 
	\begin{center}
		\includegraphics[width=60mm]{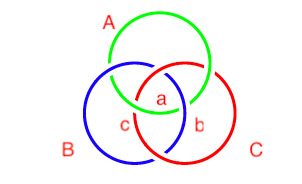}
	\end{center}
\end{figure}

From the diagram of the Borromean rings, and using the algebraic interpetration of the fundamental quandle, it is clear its  link quandle   $\Gamma$ is generated by 
$\{A,B,C,a,b,c\}$ and the generators satisfy the following relations:
\begin{equation}
	\begin{array}{lll}
		a*C=A;\, & b*A=B;\, & c*B=C;\\
		a*c=A;\, & b*a=B; \, &c*b=C;\\
	\end{array}
\end{equation}
As we can see from this relations, generators $A,B,C$ can be written in terms of the generators $a,b,c$ so they can be erased from the presentation.  So the relations turn into:
\begin{equation}
	\begin{array}{l}
		a*(c*b)=a*c;  b*(a*c)=b*a;\, c*(b*a)=c*b\\
		a*c=a*c;\, b*a=b*a; \, c*b=c*b.\\
	\end{array}
\end{equation}
Since  the relations $ a*c=a*c;\,  b*a=b*a; \, c*b=c*b$ are redundant they can be erased from the representation, so we have shown:

\begin{lemma} \label{l:pre}
	The fundamental quandle of the Borromen rings has the following  presentation:
	$$
	\langle  a,b,c: a*(c*b)=a*c;  b*(a*c)=b*a;\, c*(b*a)=c*b\rangle. 
	$$
\end{lemma}

Before we continue observe we have removed relations and generators from the initial representations in the same way we do it 
in group theory by means of Tietze transformations, we are not going to it here but it is possible to show that Tietze moves remain valid when we work with quandle presentations. 

\begin{lemma} \label{l:cor} Up to conjugacy there are exactly two injective colorations of the Borromean rings in the parabolic quandle $Par$, see example \ref{e:par}. 
\end{lemma}

{\it Proof.-} Let $\Gamma$ be the presentation of the fundamental quandle of the Borromean rings given by Lemma \ref{l:pre}.
Let $f: \Gamma\rightarrow \mathcal{P}$ be an injective quandle morphism, then $f$
satisfies:
\[
\begin{array}{ll}
	f(a)=(f(a)*(f(c)*f(b)))*^{-1}f(c)\\
	f(b)=( f( b)*(f(a)*f(c)))*^{-1}f(a)\\
	f(c)=( f( c)*(f(b)*f(a)))*^{-1}f(b).
\end{array}
\]
Since the coloration is injective and after conjugation we can assume that:
\[
\begin{array}{l}
	f(a)=[1,0],\, f(b)=[0,t],\,f(c)=[x,y];\\
\end{array}
\]
Combining the previous conditions we get: 
\[
\begin{array}{l}
	\,[1,0]=[1-x^2t^2+x^3t^2y, t^2x(t^2x-2y-t^2x^2y+xy^2)];\\
	\,[0,t]=[ty(-y+xy^2+2x-x^2y),t(y^4+1+y^2-xy^3)]; \\
	\,[x,y]=[x(1+t^2+t^4)-t^2y,y+t^4 x-t^2 y]. \\
\end{array}
\]
After some manipulations, one gets such restrictions are equivalent to solve the following system of polynomial equations:
\begin{equation} \label{e:sp}
	\begin{array}{l}
		2-x^2t^2+x^3t^2y=0\\
		t^2x-2y-t^2x^2y+xy^2=0;\\
		-y+xy^2+2x-x^2y=0\\
		y^4+2+y^2-xy^3=0; \\
		(2+t^2+t^4)x-t^2y=0\\ 
		t^4 x+(2-t^2)y =0\\
	\end{array}
\end{equation}
It is easy to see in equation \ref{e:sp}, that the last two equations form a homogeneous system of linear equations on $x,y$, which has non-trivial solutions iff the coefficient determinant is zero. A straightforward computation shows such determinant is given by $4+t^4$.  So take  $t^4=-4$, now if we substitute  this value into the last equation of  
\ref{e:sp} we get  $x=4^{-1}(2-t^2)y$,  using this relation in the equation $-y+xy^2+2x-x^2y=0$ we obtain $y^2=x^2$. After some manipulations, we conclude that  \ref{e:sp} has exactly 8 different solutions which are given by:
\[
\begin{array}{l}
	t=1+i,1-i,-1+i,-1-i\\
	y=\pm t\\
	x=4^{-1}(2-t^2)y
\end{array}
\]
Passing to $Par$ we deduce $f$ is one of the colorations:
\[
\begin{array}{lll}
	f_{1}(a)=[1,0], & f_{1}(b)=[0,1+i], & f_1(c)=[1,1+i] \\
	f_{2}(a)=[1,0],& f_{2}(b)=[0,1-i],& f_1(c)=[1,1-i] \\
\end{array}
\]
$\square$

Before we finish this section let us compute a presentation of the link group for the Borromean rings.

\begin{lemma}
	
	The fundamental group of the Borromen rings has the following  presentation:
	$$
	\langle a,b,c\vert [[c^{-1},b],a]=[[a^{-1},c],b]=[[b^{-1},a],c]=e\rangle;
	$$
\end{lemma} 

{\it Proof.-} 
From the diagram of the Borromean rings and Wirtinger Theorem we get the fundamental group of the borromean group is generate by the arcs
$\{A,B,C,a,b,c\}$ which satisfy:
\begin{equation}
	\begin{array}{lll}
		CaC^{-1}=A;\, & AbA^{-1}=B;\, & BcB^{-1}=C;\\
		cac^{-1}=A;\, & aba^{-1}=B; \, &bcb^{-1 }=C;\\
	\end{array}
\end{equation}
As we can see from this relations, generators $A,B,C$ can be written in terms of the generators $a,b,c$ so they can be erased from the presentation.  So the relations turn into:
\begin{equation}
	\begin{array}{l}
		bcb^{-1 }abc^{-1}b^{-1 }=cac^{-1};\, cac^{-1}bca^{-1}c^{-1}= aba^{-1};\, aba^{-1}cab^{-1}a^{-1}=bcb^{-1 };\\
		cac^{-1}=cac^{-1};\, aba^{-1}=aba^{-1}; \, bcb^{-1 }=bcb^{-1 }.\\
	\end{array}
\end{equation}
Since  the relations $ a*c=a*c;\,  b*a=b*a; \, c*b=c*b$ are redundant by Tietze movements they can be erased, now after some manipulations we get the relations: 
\[
\begin{array}{l}\,
	[[c^{-1},b],a]=[c^{-1},b]a[c^{-1},b]^{-1}a^{-1}=c^{-1}bcb^{-1 }abc^{-1}b^{-1 }c a^{-1}=e; \\\, 
	[[a^{-1},c],b]=[a^{-1},c]b[a^{-1},c]^{-1}b^{-1}=a^{-1}cac^{-1}bca^{-1}c^{-1}ab^{-1}= e\,\\ \,
	[[b^{-1},a],c]=[b^{-1},a]c[b^{-1},a]^{-1}c^{-1}=b^{-1} aba^{-1}cab^{-1}a^{-1}bc^{-1}=e
\end{array}
\]
which concludes the proof.
$\square$

\section{Main results}\label{s:main}
\begin{theorem}
	Let $\mathcal{L}$  the Borromean rings and $D$ its  diagram as before. Thus , up to conjugation, there are exactly two representations  
	$\rho_1,\rho_2 : \pi_1(\Bbb{S}^3\setminus \mathcal{L}) \rightarrow {\rm PSL}(2,\Bbb{C})$,
	which are faithfull, discrete and  send meridians into parabolic elements.  Moreover, we have $\rho_1=\overline{\rho_2}$.
\end{theorem}
{\it Proof.- } By Corollaries \ref{c:mod} and \ref{c:correp} the set of repressentations which are discrete, faithfull and send meridians into Parabolic elements correspond to a subset in te set of colorations of $D$ into Par which are injective. By 
Lemma \ref{l:cor}, there are two of such colorations, say $A_1, A_2$, which satisfy $A_1=\overline A_2$. To conclude observe that image of the induced representations are:

\[
\left 
\langle 
\begin{bmatrix}
	2+i & 2i\\
	-1 & -i
\end{bmatrix}, 
\begin{bmatrix}
	1 & 2i\\
	0 & 1
\end{bmatrix}, 
\begin{bmatrix}
	1 & 0\\
	-1 & 1
\end{bmatrix}
\right 
\rangle, \,\,
\left 
\langle 
\begin{bmatrix}
	2-i & -2i\\
	-1 & i
\end{bmatrix}, 
\begin{bmatrix}
	1 & -2i\\
	0 & 1
\end{bmatrix}, 
\begin{bmatrix}
	1 & 0\\
	-1 & 1
\end{bmatrix}
\right 
\rangle \subset PSL(2,\Z[i])
\]
which are discrete and by the Wirtinger theorem isomorphic to $\pi_1(\Bbb{S}^3\setminus L)$, compare with \cite{Matsumoto}.
$\square$

Now the proof of Theorem \ref{t:ext} is straighforward so we will ommit it here.

\section*{Acknowledgements}
The authors would like to thank F. González for fruitful 
conversations. The second author, also thanks staff of the UCIM at UNAM for their kindness
and 
help during his stay. This work was partially supported by grants of project PAPPIT UNAM: 
IN110219.

\bibliographystyle{plain}

\end{document}